\documentclass[12pt]{article}
\usepackage{amssymb}
\usepackage{amsmath}
\begin{document}
\def\Drb{\nobreak\hfill\nobreak\lower 1pt\hbox{$\square $}}
\def\black{\vrule height 6pt width 6pt depth 1pt} 
\newcommand{\Aut}{\mbox{Aut }}
\newcommand{\shift}{\mbox{shift}}
\newcommand{\Ad}{\mbox{Ad }}
\newcommand{\Fun}{\mbox{Fun }}
\newcommand{\ind}{\mbox{ind }}
\newcommand{\diag}{\mbox{diag }}
\newcommand{\Ind}{\mbox{Ind }}
\newcommand{\spann}{\mbox{span }}
\newcommand{\Ob}{\mbox{Ob }}
\newcommand{\Hom}{\mbox{Hom }}
\newcommand{\pfeil}{\longrightarrow}
\newcommand{\Pfeil}{\longmapsto}
\newcommand{\proj}{\mbox{proj }}
\renewcommand{\blacksquare}{\nobreak \hfill \nobreak \black}
\newcommand{\blacksquared}{\nobreak \hfill \nobreak \black}
\newcommand{\tm}{\otimes _M}
\newcommand{\tn}{\otimes _N} 
\newcommand{\tp}{\otimes _P}
\newcommand{\tq}{\otimes _Q}
\newcommand{\ta}{\otimes _A}
\renewcommand{\th}{\otimes _H}
\newcommand{\kreuz}{\rtimes}
\newcommand{\ckreuz}{{\rtimes} _c}
\newcommand{\ockreuz}{{\rtimes} _{\bar{c}}}
\newcommand{\rest}{\! \mid \!}
\newcommand{\C}[1]{{\cal #1}}
\newcommand{\B}[1]{{\bf #1}}
\newcommand{\equi}{\Longleftrightarrow}
\newcommand{\beq}{\begin{equation}}
\newcommand{\eeq}{\end{equation}}
\newcommand{\Kl}[3]{\langle #1 \otimes #2, #3\rangle}
\newcommand{\ska}[2]{\langle #1,\, #2\rangle}
\newcommand{\skab}[2]{\Bigl\langle #1,\, #2\Bigr\rangle}
\newcommand{\skal}[2]{\langle #1,\, #2\rangle _N^l}
\newcommand{\skar}[2]{\langle #1,\, #2\rangle _N^r}
\newcommand{\ov}{\overline}
\newcommand{\newl}{\newline \noindent}
\newcommand{\tr}{\mbox{tr}}
\newcommand{\Tr}{\mbox{Tr}}
\newcommand{\zws}{{\rm II}_1}
\newcommand{\abs}{\vspace{4mm}}
\newcommand{\orho}{\bar{\rho }}
\newcommand{\osigma}{\bar{\sigma }}
\newcommand{\ophi}{\bar{\phi }}
\newcommand{\opsi}{\bar{\psi }}
\newcommand{\otau}{\bar{\tau }}
\newcommand{\ovs}[1]{\ov{\sigma (#1)}}
\newcommand{\oR}{\bar{R}}
\newcommand{\oV}{\bar{V}}
\newcommand{\oW}{\bar{W}}
\newcommand{\oc}{\bar{c}}
\newcommand{\comp}{\mathbb{C}}
\newcommand{\nat}{\mathbb{N}}
\newcommand{\reel}{\mathbb{R}}
\newcommand{\ganz}{\mathbb{Z}}
\newcommand{\rat}{\mathbb{Q}}

 \title{A Note on Bimodules and $\zws $-Subfactors}
\author{Reinhard Schaflitzel\\
        Mathematisches Institut\\
        Technische Universit\"at M\"unchen\\
        Arcisstr.21\\ 
        80333 M\"unchen\\
         GERMANY\\ email: schafl@mathematik.tu-muenchen.de }
\date {September 13, 1996}

\maketitle
\newtheorem{lem}[subsection]{Lemma}
\newtheorem{cor}[subsection]{Corollary}
\newtheorem{prop}[subsection]{Proposition}
\newtheorem{thm}[subsection]{Theorem}
\newtheorem{defi}[subsection]{Definition}
\newtheorem{quest}[subsection]{Question}
\setcounter{tocdepth}{1}
\hyphenation{bundles Yamagami}

\begin{abstract} 
A brief introduction into bimodules of
$\zws$-factors is presented. Furthermore a version of the following
result due to M. Pimsner and S. Popa is derived: 
Let $N=M_{-1}\subset M=M_0 \subset M_1 \subset M_2 \subset
 \ldots  $ denote the Jones tower of a $\zws $-factor $N\subset M$
with finite index. Then the factor obtained by the basic construction  
from the pair $N\subset M_{n-1}$ is equal to $M_{2n-1}$.
\end{abstract}

\section*{Introduction} The theory of subfactors was established
by V.F.R. Jones in his famous paper \cite{Jo}.
A. Ocneanu had the idea to use bimodules (also
called correspondences) for the theory of subfactors (see \cite{Oc1} and
\cite{Oc2}). In this note we present an easy access to the bimodules of
$\zws $-factors and to their application to the theory of $\zws $-subfactors.
Let
$$N=M_{-1}\subset M=M_0 \subset M_1 \subset M_2 \subset M_3 \subset
 \ldots  $$
be the Jones tower of a $\zws $-subfactor $N\subset M$ with finite index.
The main goal of this paper is to construct a normal isomorphism
from $M_{2n-1}$ ($n\in \nat $) onto the von Neumann algebra of the
right $N$-linear operators on $L^2(M_{n-1})$. In particular
we obtain a new version of the result due to M. Pimsner and S. Popa that the
basic construction for $N\subset M_{n-1}$ is equal to $M_{2n-1}$. The author 
thinks that his result is more convenient for applications than 
Pimsner's and Popa's. We notice that
Y. Denizeau and J.F. Havet developed a theory of correspondences
of finite index for arbitrary von Neumann algebras (see \cite{DH1} and \cite{DH2}) from
which our results follow likewise. S. Yamagami also proved some results 
about bimodules and some parts of Ocneanu's approach to subfactors
(see \cite{Y1} and \cite{Y2}).
 This note contains the results of an 
introductionary section and of a part of the appendix of the author's
Habilitationsschrift \cite{Sch2}.

\section{Bimodules of $\zws $-factors}
In order to avoid
subtleties we will assume that every von Neumann algebra
appearing in this paper acts on a separable Hilbert space.
For a $\zws $-factor $S$ let $\tr _S$ denote the unique normalized 
trace of $S$.
Let $L$, $P$, $Q$, and $S$ be factors of type $\zws $. 
Let $Q^{op}$ be the factor opposite to Q. ($Q^{op}$ is equal to
$Q$ as a complex vector space, the multiplication law is reversed, that
means $p\stackrel{\rm op} \circ q:=q\cdot p$ for $p,q\in Q$, and the
involution $*$ is the same as in $Q$.)

 \begin{defi} (i)
A (left) $P$-module is a Hilbert space $\cal H$ endowed with
a normal representation
$\lambda :P\longrightarrow {\bf L}({\cal H})$ of $P$. A right $Q$-module is
a Hilbert space ${\cal H}$
endowed with a normal representation
$\rho :Q^{op}\longrightarrow {\bf L}({\cal H})$ of $Q^{op}$.
A $(P,Q)$-bimodule $\cal H$ is a left $P$-module and a right $Q$-module
such that $\lambda (P)$ and
$\rho (Q^{op})$ commute. If we like to emphasize the Hilbert space $\C{H}$, we
write $\lambda _{\C{H}}$
for $\lambda $ and $\rho _{\C{H}}$ for $\rho $.
\newl (ii) For left $P$-modules $\C{H}_1$ and $\C{H}_2$ a continuous linear
operator \linebreak $T:\C{H}_1\pfeil \C{H}_2$ is called (left) $P$-linear, if 
$T\, \lambda _{\C{H}_1}(p)= \lambda _{\C{H}_2}(p)\, T$ holds for every
$p\in P$. For right $Q$-modules $\C{H}_1$ and $\C{H}_2$ a continuous linear
operator $T:\C{H}_1\pfeil \C{H}_2$ is called right $Q$-linear, if 
$T\, \rho _{\C{H}_1}(q)= \rho _{\C{H}_2}(q)\, T$
for every $q\in Q$.
For $(P,Q)$-bimodules $\C{H}_1$ and $\C{H}_2$ an
operator  $T:\C{H}_1\pfeil \C{H}_2$ is called $(P,Q)$-linear,  
if $T$ is left $P$- and right $Q$-linear.    \label{D11y}\end{defi}

We write $_P{\cal H}_Q$ for the $(P,Q)$-bimodule $\C{H}$.
For $\xi \in {\cal H}$, $p\in P$ and $q\in Q$ we also write $p\, .\xi $ 
in place
of $\lambda (p)\xi$ and
$\xi .\, q$ in place of $\rho (q)
\xi $. The commutant $\lambda _{\C{H}}(P)'$ is also denoted by 
$\C{L}_{P,-}(\C{H})$,
the commutant $\rho _{\C{H}}(Q)'$
by $\C{L}_{-,Q}(\C{H})$, and the set 
$\lambda _{\C{H}}(P)'\cap \rho _{\C{H}}(Q)'$ of the
continuous $(P,Q)$-linear
operators on $\C{H}$ by $\C{L}_{P,Q}(\C{H})$. 

Let $L^2(Q)$ be the Hilbert space obtained by the completion of $Q$
with respect to the inner product
$\ska{x}{y}:=tr_Q(y^*x)$ ($x,y\in Q$). We denote an element $x$ of $Q$
by $\ov{x}$, if $x$ is regarded as an element of $L^2(Q)$.
$L^2(Q)$ is a $(Q,Q)$-bimodule, where the actions are given by left and
right multiplication.  

For $(P,Q)$-bimodules we use the usual concepts of representation theory.
For instance, a $(P,Q)$-bimodule $_P{\cal H}_Q$ is called irreducible, if
and only if there is no closed proper subspace
${\cal K}\neq 0$ of $\cal H$ invariant 
under the left action of $P$ and the right
action of $Q$, or equivalently, if and only if
$\C{L}_{P,Q}(\C{H}) =\comp \B{1}$. Two
$(P,Q)$-bimodules $\cal H$ and $\cal K$ are called equivalent
(${\cal H}\cong {\cal K}$),
if there is a unitary $(P,Q)$-linear map
$U:{\cal H}\longrightarrow {\cal K}$.

\pagebreak 
\subsection{The $\otimes _Q$-tensor product \label{SAQ1}}
Let $\C{H}$ be a right $Q$-modules and $\C{K}$ a left $Q$-module. We describe
J.-L. Sauvageot's construction (\cite{Sau}) of the tensor product
$\C{H}\otimes _Q\C{K} $ for our special case of 
$\zws $-factors. (It is not difficult 
to check that the definition in \cite{Su} is a special case of the 
definition in \cite{Sau}). The Hilbert
space $\C{H} \tq \C{K}$ has similar properties
like the algebraic $Q$-tensor product $\tq $ (see \cite{Bou}).

 We consider
the elements $\eta \in \C{K}$, for which there is a continuous linear operator
$R_l(\eta ):L^2(Q) \pfeil \C{K}$ such that $R_l(\eta )\, \ov{x}
\, =x\, .\eta $ for every $x
\in Q$.
These elements are called left bounded, they form a dense subspace $D_l(\C{K})$
 of $\C{K}$. For
$\eta _1,\,\eta _2 \in D_l(\C{K})$, $\skal{\eta _1}{\eta _2}:=JR_l(\eta 
_1)^*R_l(\eta _2)J$
is a right $Q$-linear continuous operator on $L^2(Q)$ and thus belongs to $Q$
(where $J$ is the antiunitary operator on $L^2(Q)$ defined by $J \ov{q} =
 \ov{q^*}$ for $q\in Q\subset L^2(Q)$.). 
So a $Q$-valued inner product
$\skal{\cdot}{\cdot}:D_l(\C{K})\times D_l(\C{K})  \pfeil Q$ is given.

For the right $Q$-module $\C{H}$ we introduce the set 
$D_r(\C{H})$ of right-bounded elements and the operator
$R_r(\xi )$ for $\xi \in D_r(\C{H})$ in an analogous way:
\newl  $R_r(\xi ):L^2(Q) \pfeil \C{H}$ is given by
$R_r(\xi ) \, \bar{x}= \xi .\, x$ for $x\in Q$, and 
$\skar{\xi _1}{\xi _2}:=R_r(\xi _2)^*
R_r(\xi _1)$ defines a $Q$-valued inner product on $D_r(\C{H})$. 
Moreover, $D_r(\C{H})$
is dense in $\C{H}$.

We use the symbol $\odot $ for algebraic tensor products
of $\comp $-vector spaces.
The Hilbert space $\C{H}\otimes _Q\C{K}$ is defined as
the Hausdorff completion of the 
algebraic tensor product $D_r(\C{H})\odot \C{K}$ with 
respect to the inner product
\begin{equation} \label{EA11}\ska{\xi _1\odot \eta _1}{\xi _2\odot \eta _2}:=
\ska{\skar{\xi _1}{\xi _2}\, .\eta _1}{\eta  _2}.  \end{equation} 
(More precisely, $\C{H}\otimes _Q\C{K}$ is the completion of the factor
space $D_r(\C{H})\odot \C{K}/N$ with respect to the inner product 
(~\ref{EA11}), where $N$ is the subspace of \linebreak
$D_r(\C{H})\odot \C{K}$
consisting of all vectors
$\psi $, for which $\ska{\psi }{\psi }=0$.)
 
You also can consider $\C{H}\tq \C{K}$ as the Hausdorff completion of 
$\C{H}\odot D_l(\C{K})$ with respect to
$$\ska{\xi _1\odot \eta _1}{\xi _2\odot \eta _2}:=\ska{\xi _1.\, \skal{\eta _1}{\eta _2}}{\xi _2}.$$
The inner products agree on $D_r(\C{H})\odot D_l(\C{K})$ (see \cite{Sau}, 
Lemma 1.7), and one easily sees that the 
image of $D_r(\C{H})\odot D_l(\C{K})$ is dense in $\C{H}\tq \C{K}$ (for both 
the definitions of $\C{H}\tq \C{K}$).

\begin{lem} The map $\xi \in \C{H} \Pfeil \xi \otimes _Q \eta $ (resp.
$\eta \in \C{K} \Pfeil \xi \otimes _Q \eta $) is continuous for every
$\eta \in D_r(\C{K})$ (resp. $\xi \in D_l(\C{H})$).
\label{LA11} \end{lem}

We call the pair $(\xi ,\eta )\in \C{H} \times \C{K}$  admissible, if
and only if $\xi \in D_r(\C{H})$ or
$\eta \in D_l(\C{K})$.
For admissible $(\xi ,\eta )$ we denote the corresponding element of
$\C{H}\otimes _Q \C{K}$ by $\xi \tq \eta $.

$D_l(\C{K})$ is stable under the left action of $Q$ and 
$\skal{q\, .\eta _1}{\eta _2}\,=q\cdot \skal{\eta _1}{\eta _2}$
for $q\in Q$ and $\eta _1,\,\eta _2 \in D_l(\C{K})$. Furthermore
$D_r(\C{H})$ is stable under the right action  of $Q$ and
$\skar{\xi _1.\, q}{\xi _2}=\skar{\xi _1}{\xi _2} \cdot q$ for 
$\xi _1,\,\xi _2 \in D_r(\C{H})$.  It follows
$$\xi .\, q \otimes _Q \eta = \xi \otimes _Q q.\, \eta $$
for $(\xi ,\eta )$ admissible and $q\in Q$. 

\abs
We have the following Lemmata:
\begin{lem}
(i) Let $\C{H}$ and $\C{H}'$ be right $Q$-modules and $\C{K}$ and $\C{K}'$ 
left $Q$-modules. 
If $A:\C{H} \pfeil \C{H}'$ is a continuous right $Q$-linear operator and 
$B:\C{K}\pfeil \C{K}'$ a
continuous left $Q$-linear operator, then $(A\tq B) \, \xi \tq \eta  
=A\xi \tq B\eta $ for admissible 
$(\xi ,\eta )$ defines a unique continuous linear operator 
$A\tq B$ from $\C{H}\tq \C{K}$ onto
$\C{H}'\tq \C{K}'$. One obtains $(A\tq B)^*=A^*\tq B^*$.
\newl (ii) $\C{H}\tq \C{K}$ is a left $P$-module with
respect to the action $\lambda (p)\tq \B{1}\, (p\in P)$, if $\C{H}$ is a 
$(P,Q)$-bimodule.
It is a right $S$-module with
respect to the right action $\B{1}\tq \rho (s)\, (s\in S)$, 
if $\C{K}$ is a $(Q,S)$-bimodule, and a $(P,S)$-bimodule, if
both conditions are satisfied.
 \newl (iii) If $\C{H}$ and $\C{H}'$ are $(P,Q)$-bimodules and $A$ is 
$(P,Q)$-linear, then
$A\tq B$ is left $P$-linear. The corresponding result for right actions on 
$\C{K}$ and $\C{K}'$ is satisfied likewise.
\label{L1} \end{lem}
         
\B{Proof:\ }We note $A(D_r(\C{H})) \subset D_r(\C{H})$ and 
$B(D_l(\C{K}))\subset D_l(\C{K})$.
Assertion (i) without the equation concerning the $*$-operation is Lemma
2.3 in \cite{Sau}. Part (ii) and (iii) are easy consequences of part (i),
hence we will only show \newl $(A\tq B)^*=A^*\tq B^*$. 
If $\C{K}=\C{K}'$ and $B=\B{1}$, then
\begin{eqnarray*} \ska{(A\tq \B{1})(\xi \tq \eta)}{\xi '\tq \eta '}\,&=&
   \ska{(A\, \xi ).\, \skal{\eta }{\eta '}}{\xi '}\, =\\
  \ska{\xi .\, \skal{\eta }{\eta '}}{A^* \, \xi '}\,
   &=&\ska{\xi \tq \eta}{(A^*\tq \B{1})(\xi '\tq \eta ')}
\end{eqnarray*}
for $\xi \in \C{H},\, \xi '\in \C{H}',\, \eta,\,\eta '\in D_l(\C{K})$. 
So we obtain $(A\tq \B{1})^*
=A^*\tq \B{1}$. Similarly $(\B{1}\tq B)^* =\B{1}\tq B^*$ follows.
Using $A\tq B= (A\tq \B{1})\cdot (\B{1}\tq B)$,
we get the assertion. \blacksquare 

\abs
Lemma ~\ref{LA11} implies:
\begin{lem} (i) If $\C{H}$ is a right $Q$-module, $\C{K}$ an $(Q,P)$-bimodule
and $\eta \in D_r(\C{K})$, then $\xi \tq \eta \in D_r(\C{H}\tq \C{K})$ for 
$\xi \in \C{H}$.
\newl (ii)  If $\C{K}_0$ is a dense linear subspace of $D_r(\C{K})$ 
and $\C{H}_0$ a dense linear subspace
of $\C{H}$, then
$$\C{H}_0 \odot _Q \C{K}_0:=
 \hbox{span} \,\,
\{\xi \tq \eta:\, \xi \in \C{H}_0,\,\eta \in \C{K}_0\}$$
is a dense subset of $D_r(\C{H}\tq \C{K})$. 
\newl (iii) If $\C{H}$ is a $(P,Q)$-bimodule and $\C{K}$ a left $Q$-module, 
one obtains the analogous results.\label{L2} \end{lem}

\begin{lem} If $\C{H}$ is a right $Q$-module,\,$\C{K}$ a $(Q,P)$-module,
and $\C{L}$ a left $P$-module, there is a unique unitary operator
$$a =a(\C{H},\C{K},\C{L})\!:\C{H}\otimes _Q(\C{K}\otimes _P\C{L}) 
\pfeil (\C{H}\otimes _Q \C{K}) \otimes _P\C{L})$$
such that
$a \, \xi \otimes _Q (\eta\otimes _P\rho )\,  =(\xi \otimes _Q\eta )
\otimes _P\rho $
for $\xi \in D_r(\C{H}),\, \eta \in \C{K}$ and $\rho \in D_l(\C{L})$. 
If $\C{H}$ is an $(L,Q)$-bimodule and $\C{L}$ a $(P,S)$-bimodule, then $a$ is
$(L,S)$-linear.
 \label{L3} \end{lem}  

\B{Proof:\ } We apply the Lemma ~\ref{L2} and see that it
suffices to check 
\begin{equation} \ska{(\xi \otimes _Q\eta )\otimes _P\rho}
{(\xi '\otimes _Q\eta ')\otimes _P\rho '}
\,=\,  \ska{\xi \otimes _Q(\eta \otimes _P\rho )}
{\xi '\otimes _Q(\eta '\otimes _P\rho ')}      
\label{sk} \end{equation}
for $\xi ,\xi '\in D_r(\C{H}),\,\eta ,\eta '\in \C{K}$ and $\rho ,\, \rho '\in D_l(\C{L})$. We compute
\begin{eqnarray*} \ska{(\xi \otimes _Q\eta )\otimes _P\rho }{(\xi '\otimes _Q
\eta ')\otimes _P\rho '}\, &=&
   \ska{(\xi \otimes _Q\eta ).\, \ska{\rho }{\rho '}_P^l}{\xi '\otimes _Q\eta '}
 \,= \\
            &=&\, \ska{\ska{\xi }{\xi '}_Q^r\, .\eta .\, 
           \ska{\rho }{\rho '}_P^l}
{\eta '}\end{eqnarray*}
and get the same result, if we carry out similar steps for the right 
inner product in (~\ref{sk}). \blacksquare 

\abs
For every $(P,Q)$-bimodule $_P{\cal H}_Q$ there is the conjugate bimodule
$_Q\bar{\cal H}_P$, where
 $\bar{\cal H}$ is equal to $\cal H$ as a real vector space, but the inner 
product and the scalar
multiplication are conjugate. The left action $\bar {\lambda }$ of $Q$ and 
the right action $\bar {\rho }$
of $P$ on $\bar{\cal H}$ are given by
$\bar{\lambda }(q)=\rho (q^*)$ for $q\in Q$ and $\bar{\rho }(p)=\lambda (p^*)$
for $p\in P$.
Obviously, ${\cal H}\cong {\cal K}$ entails $\bar{\cal H}\cong \bar{\cal K}$.

The $(Q,Q)$-bimodule $L^2(Q)$ is a unit in the following meaning:
\noindent For a $(Q,S)$-bimodule $\C{H}$ there is a unitary $(Q,S)$-linear map 
$l_{\C{H}}$ from \newl
$L^2(Q) \tq \C{H}$ onto $\C{H}$ determined by
\begin{equation} l_{\C{H}}\,\, \ov{q} \tq \xi = q\, .\xi \qquad \hbox{for $\xi \in \C{H}$ and 
     $q\in Q$.} \label{lH} \end{equation}
The existence of $l_\C{H}$ allows us to identify the bimodules
$L^2(Q) \tq \C{H}$ and $\C{H}$.
Similarly for a $(P,Q)$-bimodule $\C{H}$ there is a unitary $(P,Q)$-linear
map $r_{\C{H}}$ from $\C{H} \tq L^2(Q)$ onto $\C{H}$ determined by
\begin{equation} r_{\C{H}} \,\,
\xi \tq \ov{q} =\xi .\, q \qquad \text{for $\xi \in \C{H}$
and $q\in Q$.} \label{rH} \end{equation}

\begin{defi} A $(Q,P)$-bimodule $\C{H}$ is called regular, if 
$D(\C{H}):=D_l(\C{H})\cap D_r(\C{H})$ 
is dense in $\C{H}$. \label{DA1}\end{defi}

A subbimodule of a regular bimodule, a finite direct sum of regular 
bimodules, the 
conjugate of a regular bimodule
and the tensor product of regular bimodules are regular. 
(The proof is easy, Lemma ~\ref{L2} is used for the tensor product.)

\noindent Let $Q_i, i=0,1\ldots ,n,$ be factors of type $\zws $, and 
$_{Q_{i-1}}
\C{H}_{i\,\, {Q_i}},\,i=1,\ldots ,n$
be regular bimodules. There are different possibilities
to put the the brackets in
 $\C{H}_1 \otimes _{Q_1}\C{H}_2
\otimes _{Q_2} \cdots \otimes _{Q_{n-1}} \C{H}_n$. From Lemma ~\ref{L3}
we conclude that the $(Q_0,Q_n)$-bimodules corresponding to different
choices of the brackets are equivalent in a canonical way
and we may identify these $(Q_0,Q_n)$-bimodules.
If $D_i$ is a dense linear subset of
$D(\C{H}_i)$ for $i=1,\ldots ,n$,
then 
$$\displaylines{
D_1\odot _{Q_1}D_2\cdots \odot _{Q_{n-1}}D_n:= \hfill\cr
\hfill \spann \{\xi _1 \otimes _{Q_1}\xi _2
\otimes _{Q_2} \cdots \otimes _{Q_{n-1}} \xi _n: 
\,\xi\in D_i, i=1,\ldots ,n \}    }$$
is dense in $\C{H}_1 \otimes _{Q_1}\C{H}_2
\otimes _{Q_2} \cdots \otimes _{Q_{n-1}} \C{H}_n$. 

\section{Bimodules and the Jones Tower}
Let $N\subset M$ be a subfactor of a $\zws $-factor $M$ with finite 
Jones' index
$\beta :=[M:N]$. $N$ and $M$ act on $L^2(M)$ by left multiplication.
The von Neumann algebra $M_1:=\C{L}_{-,N}(L^2(M))$ is a $\zws $-factor
containining $M$ and is called the basic construction for $N\subset M$.
The orthogonal projection $e_0$ from $L^2(M)$ onto 
$L^2(N)\subset L^2(M)$ is called the Jones projection. It is known that
$M_1$ is  generated by $M$ and $e_0$ as a $*$-algebra.
$e_0$ maps $M\subset L^2(M)$ onto $N\subset L^2(M)$, the restriction 
$$E_0=e_0\rest M:M \pfeil N$$ of $e_0$ is a normal faithful conditional expectation
from $M$ onto $N$. For $m\in M$ $E_0(m)$ is the unique element of $N$ 
satisfying 
\[
\tr_M (E_0(m)\, n)=\tr_M (m\, n) \qquad \text{for every $n\in N$}.\]
$E_0$ is called the conditional expectation from $M$ to $N$ corresponding
to the trace $\tr _M$.

We obtain  $[M_1:M] =[M:N]
< \infty $, so we are able to repeat the basic construction 
infinitely many times and get the so called Jones
tower
$$N=M_{-1}\subset M=M_0 \subset M_1 \subset M_2 \subset M_3 \subset
 \ldots  ,$$
where $M_{k+1}$ is the basic construction for the subfactor
$M_{k-1} \subset M_k$. 
For $k\in \nat \cup \{0\}$ let
$e_k \in M_{k+1}$ denote the orthogonal
projection from $L^2(M_k)$ onto $L^2(M_{k-1})$ and $E_k:M_k \pfeil M_{k-1}$ the corresponding 
conditional expectation. We have the following relations:
\begin{eqnarray} e_k\, e_{k\pm 1}\,e_k \,&=& \beta ^{-1}e_k \qquad \hbox{and}
                     \label{E1}            \\
                  e_k\,e_l\,&=& e_l\,e_k                     \qquad \hbox{for $|k-l|\geq 2$.}\label{E2}
\end{eqnarray}
Moreover we have
\begin{eqnarray} e_k\,x\,e_k\,&=&E_k(x)e_k \qquad \hbox{for $x\in M_k$\ \ \
 and} \label{E3} \\ 
                   E_k(e_{k-1})&=& \beta ^{-1}\B{1}.  \label{E4}\end{eqnarray}

The trace $\tr _{M_{k+1}}$ on $M_{k+1}$ satisfies the Markov property
\begin{equation} \beta \,\, \tr _{M_{k+1}}(xe_k)=\tr _{M_{k+1}}(x)  \qquad 
\hbox{for $x\in M_k$.} \label{markov}  \end{equation}

$L^2(M_k)$ ($k\in \nat \cup \{0\}$) is a regular $(M,M)$-bimodule
(as well as a regular $(N,N)$-, $(N,M)$- and 
$(M,N)$-bimodule), as
$$ M_k \subset D_l(L^2(M_k)) \cap D_r(L^2(M_k))$$
shows.

\begin{lem}    
There is a unique unitary
$(M,M)$-linear operator $U_2$ from 
\linebreak $L^2(M)\otimes _N L^2(M)$ onto $L^2(M_1)$ such that
$U_2 \,\ov{m_1}\otimes _N\ov{m_2}= \beta ^{1/2}\, 
\ov{m_1e_0m_2}$ for $m_1,\,m_2\in M$. \label{L4}\end{lem}

Lemma ~\ref{L4} is a Hilbert space version of the following result 
(see \cite{GHJ},
Corollary 3.6.5): There 
is a canonical isomorphism (of algebraic $(M,M)$-bimodules) from the algebraic 
$N$-tensor product
$M\tn M$ onto $M_1$ given by $m_1\tn m_2 \Pfeil m_1e_0m_2$.

\abs
\B{Proof: }
According to Lemma ~\ref{L2} it suffices to verify 
\begin{equation}\ska{\ov{m_1}\otimes _N \ov{m_2}}{\ov{l_1}\otimes _N 
\ov{l_2}}\,=
\beta \, \tr_{M_2}((m_1e_0m_2)(l_2^*e_0l_1^*))
\label{st} \end{equation}
for $l_1,\,l_2,\,m_1,\,m_2\in M$ (observe the equations (~\ref{E3}) and 
(~\ref{markov})).
The right side of equation (~\ref{st}) is equal to
$$\beta \,\tr _{M_1}(l_1^*m_1E_0(m_2l_2^*)e_0)\,  =\,
\tr _{M}(l_1^*m_1E_0(m_2l_2^*)) \,=\,
 \ska{\ov{m_1}.\, E_0(m_2l_2^*)\, }{\ov{l_1}}.$$
So the Lemma is proved, if we show $\skal{\ov{m_2}}{\ov{l_2}}\,=E_0(m_2l_2^*)$.
 We get
$$\ska{R_l(\ov{m_2})^*\, \ov{x}}{\ov{n}}\,=\,\tr_{M}(x\, (nm_2)^*)
 \,=\, \tr_{M}(x\, m_2^*\, n^*)\,=\,\tr_{M}(E_0(xm_2^*)\, n^*)$$
for all $n\in N$ and $x\in M$.
It follows $R_l(\ov{m_2})^*\, \ov{x} =\ov{E_0(xm_2^*)}\in L^2(N)$ and
$$\skal{\ov{m_2}}{\ov{l_2}}.\, \ov{n}\,=\, J\,  R_l(\ov{m_2})^* \,\ov{n^*l_2}\,=J
\,\, \ov{E_0(n^*l_2m_2^*)}\,
=E_0(m_2\, l_2^*)\, .\ov{n}  $$   
for every $n\in N$. \blacksquare

\abs
We define unitary $(M,M)$-linear operators
$$U_k= U(k;N,M):\!L^2(M)^{\otimes _N^k}=
L^2(M)\tn \ldots \tn L^2(M)
\pfeil L^2(M_{k-1})$$
for $k\in \nat$ by induction. (Lemma ~\ref{L3} shows that it
is not necessary to use brackets in $L^2(M)\tn \ldots \tn L^2(M)$.)
Set $U_1:=id_{L^2(M)}$ and let $U_2$ be as in Lemma ~\ref{L4}.

For $k\geq 2$ let $V_k=U(k-1;M,M_1):\!L^2(M_1)^{\otimes _M^{k-1}}\pfeil 
L^2(M_{k-1})$.
Observing that the $(M,M)$-bimodules $L^2(M)$ and $L^2(M)\tm L^2(M)$ 
are equivalent in a canonical way according to (~\ref{lH}) (or (~\ref{rH})),
we use the following identification of $(M,M)$-bimodules:
\refstepcounter{equation} \label{E11Uk}
$$\displaylines{ L^2(M)^{\otimes _N^k}=\hfill  \cr
               L^2(M)\tn \left( L^2(M)\otimes_M L^2(M)\right)
                 \tn \cdots \tn  \left( L^2(M)\otimes_M L^2(M)\right) \tn L^2(M) \hfill\cr
             \hfill = 
\left( L^2(M)\otimes_N L^2(M)\right)^{\otimes _M^{k-1}}.
     \qquad (~\ref{E11Uk})  }$$
We put
$$U_k:=V_k\circ U_2^{\otimes _M^{k-1}}.$$
(Concerning the definition of the operator $U_2^{\otimes _M^{k-1}}$ see
Lemma ~\ref{L1}).
Using the notation
\begin{equation} e_{n,m}: =\left\{ \begin{array}{ll}
                    e_n\cdot e_{n+1}\cdot \ldots \cdot e_{m-1}\cdot e_m & 
           \mbox{for $n<m$,} \\
                    e_n                                                 &  
           \mbox{for $n=m$,} \\
                    e_n\cdot e_{n-1}\cdot \ldots \cdot e_{m+1}\cdot e_m &  
            \mbox{for $n>m$,}
                    \end{array}  \right. \label{enm} \end{equation}
we get
\begin{eqnarray} \lefteqn{U_k \, \ov{x_1}\tn \cdots \tn \ov{x_k} =}
      \nonumber  \\
            &=&\beta ^{k(k-1)/4}\, \ov{ x_1e_0x_2e_{1,0}x_3e_{2,0}x_4 \ldots
            x_{k-1}e_{k-2,0}x_k} \nonumber  \\
  &=&\beta ^{k(k-1)/4}\, \ov{x_1e_{0,k-2}x_2e_{0,k-3}x_3 \ldots x_{k-2}e_{0,1} 
                          x_{k-1}e_0x_k} \label{erl} \end{eqnarray}
for $x_1,\ldots,x_k \in M$.
(The first equation follows, if we identify $\ov{x_i}\in L^2(M)$ with
 $\ov{x_i} \otimes _M \ov{\B{1}}\in L^2(M) \otimes _M L^2(M)$ in 
(~\ref{E11Uk}), and so does the second
equation, if we put $\ov{x_i}=\ov{\B{1}}\otimes _M\ov{x_i}$.)
Observe that the considerations after Definition ~\ref{DA1}
imply that $\{\ov{x_1}\tn \cdots \tn \ov{x_k} :\, x_1,\ldots ,x_k \in M\}$
is dense in $L^2(M)^{\otimes _N^k}$.

\abs Identifying $L^2(M_1)$ and $L^2(M)\tn L^2(M)$, we may regard the 
orthogonal projection $e_1$ from $L^2(M_1)$ onto $L^2(M)\subset L^2(M_1)$
as a continuous linear operator of $\B{L}(L^2(M)\tn L^2(M))$ and denote
it by $F_1$.

\begin{thm} For $n\in \nat $ there is a normal isomorphism 
$J_n$ from \newl $M_{2n-1}=\C{L}_{-,M_{2n-3}}(L^2(M_{2n-2}))$ onto
$\C{L}_{-,N}(L^2(M)^{\otimes _N^n})$
satisfying the following properties:
$$\displaylines{J_n(m) =\lambda (m) \qquad \hbox{for $m\in M$,} \hfill \cr
  J_n(e_{2k}) \, \ov{x_1}\tn \ov{x_2} \ldots \tn \ov{x_n} \,=
                 \ov{x_1}\tn \ldots \tn \ov{E_0(x_{k+1})}\tn \ldots \tn 
                 \ov{x_n}\hfill \cr
  \hfill \hbox{for $x_1,\ldots ,x_n\in M$ and $k=0,1,\ldots ,n-1$,}  \cr
   J_n(e_{2k+1}) \, \ov{x_1}\tn \ov{x_2} \ldots \tn \ov{x_n} \,=
                 \ov{x_1}\tn \ldots \tn F_1(\ov{x_{k+1}}\tn \ov{x_{k+2}})
           \tn \ldots  \tn \ov{x_n}\hfill \cr
  \hfill \hbox{for $x_1,\ldots ,x_n\in M$ and $k=0,1,\ldots ,n-2$ ($n\geq
2$).}  }$$ 
\label{Jn}         \end{thm}

As $L^2(M)^{\otimes _N^n}$ and $L^2(M_{n-1})$ can be identified, Theorem
~\ref{Jn} especially states that the basic construction
$\C{L}_{-,N}(L^2(M_{n-1}))$ for $N\subset M_{n-1}$ is the same as $M_{2n-1}$.
In \cite{PP2} M. Pimsner and S. Popa proved a version of this result,
which is somewhat weaker than Theorem ~\ref{Jn}. 

Before we prove Theorem ~\ref{Jn}, we note some useful consequences.

Using $L^2(M)=L^2(M)\tm L^2(M)$ n times, we can identify $L^2(M)^{\tn ^n}$ with
\begin{equation} L^2(M)\tm L^2(M)\tn L^2(M)\tm L^2(M)\tn \ldots \tn
                L^2(M)\tm L^2(M). \label{ten} \end{equation} 
Let
\[\epsilon (k)= 
\begin{cases}0  &\text{for $k$ even,} \\
             -1 &\text{for $k$ odd.} \end{cases} \]
For $0\leq k < l \leq 2n$ let $\C{H}_k^l$ denote the tensor product
$$ L^2(M) \otimes _{M_{\epsilon (k)}}L^2(M) \otimes _{M_{\epsilon (k+1)}} 
\ldots \otimes _{M_{\epsilon (l-2)}} 
L^2(M),$$
which consists of the $(k+1)$.th, $(k+2)$.th, $\ldots $ and $l$.th factor in the tensor product (~\ref{ten}),
 moreover let $\C{H}_0^0:= 
\C{H}_{2n}^{2n}:=L^2(N)$.

\begin{cor} $J_n:M_{2n-1}\pfeil \C{L}_{-,N}((L^2(M)\tm L^2(M))^{\otimes _N^n})$ satisfies the following
relations:
\newl (i) $J_n(M_k) =\C{L}_{-,M_{\epsilon (k)}}(\C{H}_0^{k+1}) 
\otimes _{M_{\epsilon (k)}} \comp \B{1}_{\C{H}_{k+1}^{2n}}$
          for $-1\leq k\leq 2n-1$,
\vspace{2pt} \newl (ii) 
$J_n(M_k'\cap M_l) =\comp \B{1}_{\C{H}_0^{k+1}} \otimes _{M_{\epsilon (k)}}
            \C{L}_{M_{\epsilon (k)},M_{\epsilon (l)}}(\C{H}_{k+1}^{l+1})  
           \otimes _{M_{\epsilon (l)}} \comp
\B{1}_{\C{H}_{l+1}^{2n}}$ for \newl $-1\leq k
  <l\leq 2n-1$.
\label{Jnc}               \end{cor}

\B{Proof\ } (i) The case $k=2r-1 \,(r\in \nat ,\, r\leq n)$ follows from
$J_n(M_k)=J_r(M_k) \tn \comp \B{1}_{\C{H}_{k+1}^{2n}}$ and
$J_r(M_k)=\C{L}_{-,N}(\C{H}_0^{k+1})$ (according to Theorem ~\ref{Jn}).
In the case $k=2r, \,( r\in \nat ,\, r<n)$ we know that $J_n(M_{k-1}) \cup J_n(e_{k-1})$
and consequently $J_n(M_k)$ is contained in 
$\C{L}_{-,M}(\C{H}_0^{k+1})\tm \comp \B{1}_{\C{H}_{k+1}^{2n}}$. 
The remaining 
inclusion is a consequence of
\begin{eqnarray*} [J_n(M_{k+1}):J_n(M_k)]\,& =& [M:N] = 
     [\C{L}_{-,N}(L^2(M)^{\otimes _N^{r+1}}):
                 \C{L}_{-,M}(L^2(M)^{\otimes _N^{r+1}})] \\
                   &=&[J_n(M_k):
                 \C{L}_{-,M}(\C{H}_0^{k+1})\tm \comp \B{1}_{\C{H}_{k+1}^{2n}}] .
\end{eqnarray*}
(ii) follows from (i) and from the formula
$$(\C{L}_{-,Q}(\C{H})\tq \comp \B{1}_{\C{K}})' = \comp \B{1}_{\C{H}}\tq 
\C{L}_{Q,-}(\C{K}) $$
for a $\zws $-factor $Q$, a right $Q$-module $\C{H}$ and a left
$Q$-module $\C{K}$.
\blacksquare

\abs
Corollary ~\ref{Jnc} gives a useful presentation of the standard invariant
\begin{equation} \begin{array}{rcccccccc}
\comp \B{1}=N'\cap N &\subset &N' \cap M & \subset &N'\cap M_1&\subset & N'\cap M_2 &\subset & 
\ldots \\
 & &\cup & &\cup & &\cup & & \\
\comp \B{1}& =&M'\cap M & \subset &M' \cap M_1& \subset &M' \cap M_2 &\subset
&\ldots 
\end{array}  \label{stand} \end{equation}
of the subfactor $N\subset M$.

For a projection $p$ in $L^2(M)^{\tn ^k}$ the Hilbert space
$pL^2(M)^{\tn ^k}$ is an $(N,N)$-subbimodule of
$L^2(M)^{\tn ^k}$, if and only if $p\in J_k( N'\cap M_{2k-1})$, 
and $pL^2(M)^{\tn ^k}$ is an
 irreducible $(N,N)$-bimodule, if and only if $J_k^{-1}(p)$ is a 
minimal projection of $N'\cap M_{2k-1}$.
Now let $f$ and $g$ denote two minimal projections of $N'\cap M_{2k-1}$.
One easily sees that the $(N,N)$-bimodules $J_k(f)L^2(M)^{\tn k}$ and
$J_k(g)L^2(M)^{\tn k}$ are equivalent, if and only if
$f$ and $g$ belong to the same simple direct summand of $N'\cap M_{2k-1}$.
So we get a bijective correspondence between the equivalence classes of the
irreducible $(N,N)$-bimodules
contained in $L^2(M)^{\tn ^k}$ and the simple direct summands of 
$N'\cap M_{2k-1}$.

If one considers the principal graph of the subfactor $N\subset M$, which 
contains the information about the tower of the finite dimensional
von Neumann algebras in the upper line of (~\ref{stand}), then one usually
identifies the simple direct summands of $N'\cap M_{2k-1}$ with simple
direct summands of $N'\cap M_{2k+1}$ such that the following holds:
If $f$ is a 
minimal projection
of  a simple direct summand of $N'\cap M_{2k-1}$, then 
$fe_{2k}$ is a minimal projection of the corresponding direct summand of
$N'\cap M_{2k+1}$ (compare \cite{GHJ}, Section 4.6).

Using the description of the algebras $N'\cap M_{2k-1}$ and $N'\cap M_{2k+1}$ 
in Corollary ~\ref{Jnc}, we see that corresponding simple direct
summands of $N'\cap M_{2k-1}$ and $N'\cap M_{2k+1}$ belong to the same 
equivalence class of irreducible $(N,N)$-bimodules: Let $f$ be a minimal
projection of $N'\cap M_{2k-1}$. 
The $(N,N)$-bimodule $J_k(f)L^2(M)^{\tn ^k}$ is equivalent to
the $(N,N)$-bimodule
\begin{eqnarray*} J_k(f)L^2(M)^{\tn ^k}\tn L^2(N) &=& J_k(f)L^2(M)^{\tn ^k}\tn e_0 L^2(M) \\
   &=&  J_{k+1}(fe_{2k})L^2(M)^{\tn ^{k+1}} \end{eqnarray*}
by Theorem ~\ref{Jn}. 

In a similar way the simple direct summands of $N'\cap M_{2k}$ ($k\in \nat
\cup \{0\}$) correspond to equivalence classes of irreducible
$(N,M)$-bimodules. As above, given a minimal projection $f$ of 
$N'\cap M_{2k}$, the irreducible $(N,M)$-bimodules 
$J_{k+1}(f)(L^2(M)^{\tn ^{k+1}})$ and
$J_{k+2}(fe_{2k+1})(L^2(M)^{\tn ^{k+2}})$
are equivalent. Obviously one obtains an analogous connection between the
algebras $M'\cap M_{2k}$ and the $(M,M)$-bimodules as well as between
the
algebras $M'\cap M_{2k+1}$ and the $(M,N)$-bimodules. 

\abs
For the proof of Theorem ~\ref{Jn} we need some preparations.
\subsection {Some observations about bimodules with different 
factors acting from right \label{Bi}}
(1) Let $_Q\C{H}_L$ and $_Q\C{K}_S$ be two bimodules.
We assume that  there is an isomorphism 
$I$ from $\C{L}_{-,L}(\C{H})$ onto $\C{L}_{-,S}(\C{K})$ such that 
$I(\lambda _{\C{H}}(q))=\lambda _{\C{K}}(q)$ for every
$q\in Q$. Then $\C{H}$ is also an $\C{L}_{-,S}(\C{K})$-module, where the action
is given by $I^{-1}$. Let us suppose that the $\C{L}_{-,S}(\C{K})$-module
$\C{K}$ is equivalent to a submodule of $\C{H}$. 

A linear isometry $W:\C{K}\pfeil \C{H}$ is said to be associated with $I$, 
if $W$ is
$\C{L}_{-,S}(\C{K})$-linear. Now let us fix a linear isometry $W$ associated 
with
$I$. Let $p$ be the projection $W\cdot W^*\in \C{L}_{-,L}(\C{H})'=
\rho (L)$. We have
\begin{equation} W\, a\, W^* =p\, I^{-1}(a) \label{EA12} \end{equation}
for every $a\in \C{L}_{-,S}(\C{K})$. By restricting $W$ to its image
$p\C{H}$, we obtain a unitary operator $\tilde{W}:\C{K}\pfeil p\C{H}$.
If we endow $p\C{H}$ with the right action 
$\rho _{p\C{H} }$ of $S$ defined by $\rho _{p\C{H} }(s):=\tilde{W}\cdot \rho (s)
\cdot \tilde{W}^*$ for $s\in S$, then $p\C{H}$ is a $(Q,S)$-bimodule, 
equivalent to $\C{K}$. Using (~\ref{EA12}) we get
\begin{equation}
\C{L}_{-,S}(p\C{H})=\tilde{W}\,\C{L}_{-,S}(\C{K})\,\tilde{W}^* = W\,I(\C{L}_{-,L}(\C{H}))\,W^*
  =p\, \C{L}_{-,L}(\C{H}), \end{equation}
which implies
\begin{equation} \rho _{p\C{H}}(S)= p\, \rho _{\C{H}}(L)\, p. 
\label{rho} \end{equation}

\abs
\noindent (2) Additionally, let $_P\C{L}_Q$ be a bimodule. 
Then $\C{L} \tq \C{H}$ is a 
$(P,L)$-bimodule and $\C{L}\tq \C{K}$ is a $(P,S)$-bimodule. 
Starting with the isomorphism $I$ and
the linear isometry $W$, we will define an isomorphism 
$$\B{1}\tq I:\C{L}_{-,L}(\C{L}\tq \C{H}) \pfeil \C{L}_{-,S}(\C{L}\tq \C{K}) $$
with the following properties:
\begin{itemize}
\item[(a)] $(\B{1}\tq I)\, (\lambda _{\C{L}\tq \C{H}}(x))\, =
\lambda _{\C{L}\tq \C{K}}(x)$ for every $x\in P$ and
\item[(b)] The $\C{L}_{-,S}(\C{L}\tq \C{K})$-module $\C{L}\tq \C{K}$ is 
equivalent to a submodule of the $\C{L}_{-,S}(\C{L}\tq \C{K})$-module 
$\C{L}\tq \C{H}$, and $\B{1}\tq W$ is a linear isometry associated with
$\B{1}\tq I$. 
\end{itemize}
$\B{1}\tq p$ is a projection of $\rho _{\C{L}\tq \C{H}}(L)$, hence 
$x\Pfeil (\B{1}\tq p)\, x$ defines an isomorphism $J_1$ from 
$\rho _{\C{L}\tq \C{H}}(L)'= \C{L}_{-,L}(\C{L}\tq \C{H})$ onto the commutant
\linebreak
$((\B{1}\tq p)\, \rho _{\C{L}\tq \C{H}}(L) \,(\B{1}\tq p))'$
\ (in $(\B{1}\tq p)\, \C{L}\tq \C{H}=\, \C{L}\tq p\C{H}$ ).

In (1) $p\C{H}$ was endowed with a $(Q,S)$-bimodule structure, hence
$\C{L}\tq p\C{H}$ is a $(P,S)$-bimodule. By applying equation
(~\ref{rho}) we get
$$\rho _{\C{L}\tq p\C{H}}(S)= \B{1}\tq \rho_{p\C{H}}(S)=\B{1}\tq p
\rho _{\C{H}}(L)p  =
(\B{1}\tq p) \rho_{\C{L}\tq \C{H}}(L) (\B{1}\tq p)$$ 
with the consequence that 
$$((\B{1}\tq p)\, \rho _{\C{L}\tq \C{H}}(L) \,(\B{1}\tq p))'=
\C{L}_{-,S}(\C{L}\tq p\C{H}).$$ 
Now $x\Pfeil (\B{1}\tq \tilde{W})^*\, x\, (\B{1}\tq \tilde{W})$ 
defines an isomorphism $J_2$ from \newl 
$\C{L}_{-,S}(\C{L}\tq p\C{H})$ onto
$\C{L}_{-,S}(\C{L}\tq \C{K})$. We put $\B{1}\tq I:=J_2\circ J_1$. Property
(a) is obvious, for the proof of Property (b) we note that
$(\B{1}\tq W)\cdot (\B{1}\tq W)^* =\B{1}\tq p$ belongs to
$\rho_{\C{L}\tq \C{H}}(L)$. From
$$\displaylines{
(\B{1}\tq W)\,(\B{1}\tq I)(a)\, = (\B{1}\tq W)\, (\B{1}\tq \tilde{W})^*\,
(\B{1}\tq p)\, a\,(\B{1}\tq \tilde{W})\, = \hfill\cr
(\B{1}\tq p)\, a\, (\B{1}\tq W)\,\hspace{5mm} =a\, (\B{1}\tq p)\, (\B{1}\tq W) \,=
a  \, (\B{1}\tq p)
\hfill }$$
for every $a\in \C{L}_{-,L}(\C{L}\tq \C{H})$ we conclude that
the linear isometry $\B{1}\tq W$ is associated with $\B{1}\tq I$.

We point out that the definition of $\B{1}\tq I$ does not depend on the choice
of $W$. 

\abs
\noindent (3) We suppose the assumptions from (1). Additionally, let 
$_Q\C{M}_T$ be a bimodule and 
$J: \C{L}_{-,S}(\C{K})\pfeil \C{L}_{-,T}(\C{M})$ an isomorphism such that
$\C{M}$ is an $\C{L}_{-,T}(\C{M})$-submodule of $\C{K}$.
Considering the isomorphism 
$J\circ I:\C{L}_{-,L}(\C{H})\pfeil \C{L}_{-,T}(\C{M})$
we also may regard $\C{M}$ as an $\C{L}_{-,T}(\C{M})$-subbimodule of $\C{H}$.
If $V:\C{M} \pfeil \C{K}$ is a linear isometry associated with $J$,
then $W\circ V$ is a linear isometry associated with $J\circ I$. The proof
of this fact is easy.

\abs
\begin{lem} (i) $\beta \, e_{n,0}e_{n+1,0} =e_{n,0}e_{n+1,2}$ for $n\in \nat $.
\newl (ii) $\beta ^{k(k+1)/2} e_{n,0}e_{n+1,0}\cdot \ldots \cdot e_{n+k,0} =
            e_{n,0}e_{n+1,2}e_{n+2,4}\cdot \ldots \cdot e_{n+k,2k}$ for $0 < k\leq n $.
\newl (iii) $e_{n,0}e_{n+1,0}\cdot \ldots \cdot e_{2n-2,0}e_{2n-1,2n-2} = \beta ^{n-1}
            e_{n,0}e_{n+1,0}\cdot \ldots \cdot e_{2n-1,0}$ for $n\in \nat$.
\label{eee} \end{lem}

\B{Proof:\ } $\beta \,e_{n,0}e_{n+1,0}=\,\beta e_{n,1}e_{n+1,2}e_0e_1e_0
      =e_{n,1}e_{n+1,2}e_0 =e_{n,0}e_{n+1,2}$ shows (i). 
We get (ii) by applying (i) several times and
(iii) by applying (ii) twice (with $k=n-1$ to the right
side of the equation in (iii) and with $k=n-2$ to the left side). \blacksquare 

\subsection{The isomorphisms $J_n$ \label{Jnn}}
(1) We introduce an isomorphism 
$$I:M_3=\C{L}_{-,M_1}(L^2(M_2)) \pfeil \C{L}_{-,N}(L^2(M_1))$$
satisfying $I(m)=\lambda (m)$ for every $m\in M$. 
Let $K$ be the canonical isomorphism from $\C{L}_{-,M_1}(L^2(M_1))=
\lambda (M_1)$  onto
$\C{L}_{-,N}(L^2(M))=M_1$. $K$ is an isomorphism as in Section ~\ref{Bi}
(1) (with $Q=M$) and 
$$V:L^2(M) \pfeil L^2(M_1),\,\ov{m}\Pfeil \beta ^{1/2}\, \ov{me_0},$$
is a linear isometry associated with $K$.
(The Markov property (~\ref{markov}) of $\tr _{M_1}$ implies 
$$ \ska{\beta ^{1/2}\, \ov{me_0}}{\beta ^{1/2}\, \ov{le_0}}\, =
\beta \, \tr _{M_1}(me_0l^*)\, =\tr _M(ml^*)\, =\ska{\ov{m}}{\ov{l}}$$
for $m,l\in M$, hence $V$ is isometric. $V$ is $M$-linear and  
$V\circ e_0\, =\lambda (e_0)\circ V$ holds by equation (~\ref{E3}), 
thus $V$ is $M_1$-linear.) 

Applying Section ~\ref{Bi} (2)
with the $(M,M)$-bimodule $\C{L}=L^2(M_1)$ we get an isomorphism
$$\B{1}\otimes _M K:\C{L}_{-,M_1}(L^2(M_1)\tm L^2(M_1)) \pfeil 
\C{L}_{-,N}(L^2(M_1)\tm L^2(M)).$$
We identify $L^2(M_1)\tm L^2(M_1)$ and
$L^2(M_2)$ according to Lemma ~\ref{L4} as well as 
$L^2(M_1)\tm L^2(M)$ and $L^2(M_1)$ according to relation (~\ref{rH}).
After these
identifications, $\B{1}\tm K$ is the desired isomorphism $I$, and the linear 
isometry
$\B{1}\tm V:L^2(M_1) \pfeil  L^2(M_2)$ is associated with $I$ and
satisfies the relation
\begin{equation} (\B{1}\tm V)\, \ov{x} =\beta \,\ov{xe_{1,0}} \qquad \hbox{for } x\in M_1. 
\label{EA110} \end{equation}

\medskip\noindent (2)
Inductively  we define isomorphisms
$$J_n:M_{2n-1}=\C{L}_{-,M_{2n-3}}(L^2(M_{2n-2})) \pfeil \C{L}_{-,N}(L^2(M)^{\otimes _N^n})\qquad \hbox{($n\geq 1$)} $$
with the following properties:
\newl (a) $J_n(m)=\lambda (m)$ for every $m\in M$ and

$$\displaylines{(b) \hbox{\ \  }
W_n:L^2(M)^{\tn ^n} \pfeil L^2(M_{2n-2}),\, \, 
\ov{x_1}\tn \ov{x_2}\tn \ldots \tn \ov{x_n} \Pfeil\hfill\cr
     \hfill    \beta ^{(n-1)(2n-1)/2} \, \ov{x_1e_0x_2e_{1,0}x_3 \ldots x_{n-1}e_{n-2,0}x_ne_{n-1,0}e_{n,0}\ldots
                                    e_{2n-3,0}}\, = \cr
 \hfill \beta ^{(n-1)(2n-1)/2} \, \ov{x_1e_{0,2n-3}x_2e_{0,2n-4}x_3 \ldots x_{n-1}
                                e_{0,n-1}x_ne_{0,n-2}e_{0,n-3}\ldots
                                    e_{0,1}e_0} ,}$$
is a linear isometry associated with $J_n$ for $n\geq 2$. 

Let $J_1= id_{M_1}$ and let $W_1=id_{L^2(M)}$. 

Now for $n\geq 2$ we assume that $J_{n-1}$ is defined and that $W_{n-1}$ is associated with $J_{n-1}$. Let
$$I_n:M_{2n-1} =\C{L}_{-,M_{2n-3}}(L^2(M_{2n-2})) \pfeil \C{L}_{-,M_{2n-5}}(L^2(M_{2n-3}))$$
denote the isomorphism $I$ from Part 
(1) and $V_n$ the linear isometry $V$ from (1),
where $N$ is replaced by $M_{2n-5}$, $M$ by $M_{2n-4}$ and so on. 
We write $\hat{V}_n $ in place of $\B{1}\otimes _{M_{2n-4}}V_n$.

Observe that $R_{2n-4}:= (id_{L^2(M)}\tn U_{2n-4})\circ U_{2n-3}^*$
is an $(M,M)$-linear unitary
operator from $L^2(M_{2n-3})$ onto $L^2(M)\tn L^2(M_{2n-4})$ 
and that $$R_{2n-4}\, \beta ^{2n-3/2}\, \ov{xe_{0,2n-4}y}\, =\ov{x}\tn \ov{y}$$
holds for $x\in M$ and $y\in M_{2n-4}$. Hence it is possible to identify
the bimodules $L^2(M_{2n-3})$ and $L^2(M)\tn L^2(M_{2n-4})$.
Applying Section ~\ref{Bi} (2) we get an isomorphism $\B{1} \tn J_{n-1}$ from
$\C{L}_{-,M_{2n-5}}(L^2(M) \tn L^2(M_{2n-4}))=
 \C{L}_{-,M_{2n-5}}(L^2(M_{2n-3})) $ onto
$\C{L}_{-,N}(L^2(M)^{\otimes _N^n})$ 
and  put 
$$J_n:=(\B{1}\tn J_{n-1})\circ I_n.$$
Obviously $J_n$ satisfies Property (a).
From Section ~\ref{Bi} (3) we know that
$W_n:=\hat{V}_n\circ (\B{1}\tn W_{n-1})$ is associated with $J_n$. 
The following computation shows that $W_n$  
fulfils  Property (b) for $n\geq 3$:
$$\displaylines{(\hat{V}_n \circ  (\B{1}\tn W_{n-1}))\, 
\ov{x_1}\tn \ov{x_2}\tn \ldots \tn \ov{x_n}\,=\hfill \cr
\beta ^{(n-2)(2n-3)/2} \, \hat{V}_n \,\ov{x_1} \tn \ov{x_2e_0x_3\ldots x_{n-1}e_{n-3,0}x_ne_{n-2,0}e_{n-1,0}
            \ldots e_{2n-5,0}}  =
\hfill\cr
\beta ^{(n-1)(2n-3)/2\, +1} \,\ov{x_1e_{0,2n-4}x_2e_0x_3\ldots x_{n-1}e_{n-3,0}
x_ne_{n-2,0}e_{n-1,0}\ldots}\hfill\cr
    \hfill \ov{\ldots e_{2n-5,0}e_{2n-3,2n-4}}
=\cr
\beta ^{(n-1)(2n-1)/2} \, 
\ov{x_1e_{0,2n-4}x_2e_0x_3\ldots x_{n-1}e_{n-3,0}x_ne_{n-2,0}e_{n-1,0}
      \ldots e_{2n-5,0}e_{2n-3,0}}
      = \hfill\cr
\beta ^{(n-1)(2n-1)/2} \,
\ov{x_1e_{0,2n-3}x_2e_0x_3\ldots x_{n-1}e_{n-3,0}x_ne_{n-2,0}e_{n-1,0}\ldots e_{2n-5,0}e_{2n-4,0}}
=\hfill\cr
\beta ^{(n-1)(2n-1)/2} \,
\ov{x_1e_{0,2n-3}x_2e_{0,2n-4}x_3\ldots 
x_{n-1}e_{0,n-1}x_ne_{0,n-2}e_{0,n-3}\ldots e_0}.\hfill }$$
(We used equation (~\ref{EA110}) in line 3, Lemma ~\ref{eee} (iii) 
in line 4 and equation (~\ref{erl}) in line 6.)
\newl The case $n=2$ follows from $W_2 =\B{1} \otimes _M V$.

\subsection{Proof of Theorem 2.2}
We prove that $J_n$ satisfies the properties required in Theorem ~\ref{Jn}.
It suffices to show 
\refstepcounter{equation} \label{prv}
$$\displaylines{
 e_l\, .\,  W_n ( \ov{x_1}\tn \ov{x_2}\tn \ldots \tn \ov{x_n})\,= \hfill\cr
\,\,\,     =     \left\{ \begin{array}{ll}
     W_n\,\, \ov{x_1}\tn \ldots \ov{E_0(x_{k+1})}\tn 
     \ldots \tn \ov{x_n} &\mbox{($l=2k$ even)}\\
     W_n\, \, \ov{x_1}\tn \ldots F_1(\ov{x_{k+1}}\tn \ov{x_{k+2}})\tn 
    \ldots \tn \ov{x_n} &
     \mbox{($l=2k+1$ odd)} \end{array} \right.\hfill  (~\ref{prv}) }$$ 
for $l\in \nat \cup \{0\}, \,2n-2\geq l$ and $x_1,\ldots ,x_n\in M$.
We show (~\ref{prv}) by induction over $l$ 
(simultaneously for every $n\in \nat$ with $2n-2\geq l$):
\newl $l=0$ is a simple consequence of $e_0x_1e_0 =E_0(x_1)e_0$.

\vspace{3mm} \noindent $l=1:$ We get
\begin{eqnarray*} e_1\, .\, W_n \, \ov{x_1}\tn \ov{x_2}\tn 
                 \ldots \tn \ov{x_n} &=&
               \beta ^{(n-1)(2n-1)/2}\, e_1\, . \,
               \ov{x_1e_0x_2e_{1,0}\ldots  }\,\,\, =\\
               \beta ^{(n-1)(2n-1)/2}\,  \ov{E_1(x_1e_0x_2)e_{1,0}\ldots  }& =&
               W_n \, F_1(\ov{x_1}\tn \ov{x_2})\tn \ldots \tn \ov{x_n}.    
\end{eqnarray*}
Concerning the last $'='$ observe that there are an $m\in \nat $ and
$y_j,\, z_j,$ \newl $j=1,\ldots ,m$, such that
$E_1(x_1e_0x_2)=\sum _{j=1}^m y_j e_0 z_j$ and consequently
\newl $F_1(\ov{x_1}\tn \ov{x_2}) =\sum _{j=1}^n \ov{y_j} \tn \ov{z_j}$
holds.

\vspace{3mm}\noindent $2\leq l \leq 2n-2$:
 The left side of Equation (~\ref{prv}) is equal to
$$                  \beta ^{(n-1)(2n-1)/2} e_l\, .\,
              \ov{x_1e_{0,2n-3}x_2e_0x_3\ldots x_{n-1}e_{n-3,0}x_ne_{n-2,0}e_{n-1,0}
               \ldots e_{2n-5,0}e_{2n-4,0}} $$
(see the computation from Section ~\ref{Jnn} (2)).
For the right side we use the facts
that $W_n=\hat{V}_n\circ (\B{1}\tn W_{n-1})$ holds and that the assertion is
fulfilled for $l-2$ and $n-1$ by induction hypothesis.
For $l<2n-2$ the right side is equal to
$$\displaylines{\hat{V}_n \,\, \ov{x_1} \tn 
\Bigl( e_{l-2}\, .\,W_{n-1}\, \ov{x_2}\tn \ldots \tn \ov{x_n}\Bigr) \,=
\hfill\cr
  \beta ^{(n-2)(2n-3)/2}\, \hat{V}_n \,\ov{x_1} \tn  e_{l-2}\, . \,
 \ov{x_2e_0x_3\ldots x_{n-1}e_{n-3,0}x_ne_{n-2,0}e_{n-1,0}\ldots e_{2n-5,0}}  
\hfill = \cr
  \beta ^{(n-1)(2n-1)/2}\,  \ov{x_1e_{0,2n-3}e_{l-2}
     x_2e_0x_3\ldots x_{n-1}e_{n-3,0}x_ne_{n-2,0}e_{n-1,0}\ldots e_{2n-5,0}e_{2n-4,0}} \hfill  }$$
(The computation is similar to that in Section ~\ref{Jnn} (2).) 
So the left side 
and the right side of equation (~\ref{prv}) coincide, if
$e_le_{0,2n-3}= e_{0,2n-3}e_{l-2}$, which, however, the following computation
shows:
\begin{eqnarray*} e_le_{0,2n-3} =e_{0,l-2}e_le_{l-1}e_le_{l+1}\ldots e_{2n-3} &=& \beta ^{-1}e_{0,l-2}
                  e_{l,2n-3} = \\
                  e_{0,l-2}e_{l-1}e_{l-2}e_{l,2n-3} &=& e_{0,2n-3}e_{l-2}. 
\end{eqnarray*}
At last we deal with the case $l=2n-2$. Using the induction hypothesis and
similar arguments as before, we see that the right side of
(~\ref{prv}) is equal to
$$\displaylines{\beta ^{(n-2)(2n-3)/2}\,
  \hat{V}_n \,\ov{x_1} \tn \ov{E_{2n-4}(
    x_2e_{0,2n-5}x_3\ldots x_{n-1}e_{0,n-2}x_ne_{0,n-3}e_{0,n-4}\ldots e_0)}
\,     =\hfill \cr
        \beta ^{(n-1)(2n-3)/2 +1}\,  \ov{x_1e_{0,2n-4}E_{2n-4}(
    x_2e_{0,2n-5}x_3 \ldots } \hfill\cr
     \hfill \ov{\ldots x_{n-1}e_{0,n-2}x_ne_{0,n-3}e_{0,n-4}\ldots e_0)e_{2n-3,2n-4}}  = \cr
     \beta ^{(n-1)(2n-3)/2 +1} \, \ov{x_1e_{0,2n-4}
    x_2e_{0,2n-5}x_3\ldots } \hfill \cr 
     \hfill \ov{\ldots x_{n-1}e_{0,n-2}x_ne_{0,n-3}e_{0,n-4}\ldots  
        e_0e_{2n-4}e_{2n-3}e_{2n-4}} \, = \cr
   \beta ^{(n-1)(2n-3)/2 } \, \ov{x_1e_0
    x_2e_{1,0}x_3\ldots x_{n-1}e_{n-2,0}x_ne_{n-1,0}\ldots 
         e_{2n-4,0}e_{2n-4}} \, =\hfill\cr
   \beta ^{(n-1)(2n-3)/2 +1} \, \ov{E_{2n-2}(x_1e_0
    x_2e_{1,0}x_3\ldots x_{n-1}e_{n-2,0}x_ne_{n-1,0}\ldots}\hfill\cr
     \hfill \ov{\ldots e_{2n-4,0}e_{2n-3}e_{2n-4})}  =\, \cr
   \beta ^{(n-1)(2n-1)/2} \, e_{2n-2}\, . \,\ov{x_1e_0
    x_2e_{1,0}x_3\ldots x_{n-1}e_{n-2,0}x_ne_{n-1,0}\ldots 
     e_{2n-4,0}e_{2n-3,0}}. \hfill  }$$
(We used (~\ref{E3}) in line 3, (~\ref{E1}) in line 4, 
(\ref{E4}) in line 5, and Lemma ~\ref{eee} 
(iii) in line 6.) \blacksquare 
\pagebreak

\end{document}